\documentclass[leqno,12pt]{amsart}
\usepackage{amssymb,amscd}

\calclayout

\newtheorem{theorem}{Theorem}
\newtheorem{corollary}{Corollary}

\newtheorem{lemma}{Lemma}

\newcommand{\cF}{{\mathcal F}}
\newcommand{\cG}{{\mathcal G}}
\newcommand{\cI}{{\mathcal I}}
\newcommand{\cL}{{\mathcal L}}
\newcommand{\cM}{{\mathcal M}}

\newcommand{\cO}{{\mathcal O}}
\newcommand{\cX}{{\mathcal X}}
\newcommand{\cY}{{\mathcal Y}}

\newcommand{\mA}{{\mathbb A}}
\newcommand{\mC}{{\mathbb C}}
\newcommand{\mP}{{\mathbb P}}
\newcommand{\mZ}{{\mathbb Z}}

\newcommand{\oG}{{\overline{G}}}
\newcommand{\oP}{{\overline{P}}}
\newcommand{\oT}{{\overline{T}}}

\newcommand{\tX}{{\tilde X}}
\newcommand{\tY}{{\tilde Y}}
\newcommand{\tZ}{{\tilde Z}}

\newcommand{\ad}{\operatorname{ad}}
\newcommand{\codim}{\operatorname{codim}}
\newcommand{\diag}{\operatorname{diag}}
\newcommand{\reg}{\operatorname{reg}}

\title{Positivity in the Grothendieck group of complex flag varieties}
\author{Michel~Brion}
\address{Universit\'e de Grenoble I\\
D\'epartement de Math\'ematiques\\
Institut Fourier, UMR 5582 du CNRS\\
38402 Saint-Martin d'H\`eres Cedex, France}
\email{Michel.Brion@ujf-grenoble.fr}

\date{}
\begin{document}

\begin{abstract}
We prove a conjecture of A. S. Buch concerning the structure constants
of the Grothendieck ring of a flag variety with respect to its basis
of Schubert structure sheaves. For this, we show that the coefficients
in this basis of the structure sheaf of any subvariety with rational
singularities, have alternating signs. Equivalently, the class of the
dualizing sheaf of such a subvariety is a nonnegative combination of
classes of dualizing sheaves of Schubert varieties.
\end{abstract}

\maketitle

\section*{Introduction}

Consider a complex flag variety $X$ (see the end of this introduction
for detailed notation and conventions). The Chow ring $A^*(X)$ has an
additive basis consisting of the classes $[X_w]$ of Schubert
subvarieties. It is well known that the structure constants of
$A^*(X)$ with respect to this basis are positive, i.e.,
$$
[X_u] \cdot [X_v]=\sum_w \, a_{u,v}^w \, [X_w]
$$
for nonnegative integers $a_{u,v}^w$. A generalization to the
Grothendieck ring $K(X)$ of vector bundles (or of coherent sheaves) on
$X$ was recently formulated by A.~S.~Buch, see \cite{Bu}. The classes
$[\cO_{X_w}]$ of structure sheaves of Schubert subvarieties form an
additive basis of $K(X)$ ; this defines integers $c_{u,v}^w$ such that
$$
[\cO_{X_u}]\cdot [\cO_{X_v}]=\sum_w \, c_{u,v}^w \, [\cO_{X_w}].
$$
Setting $N(u,v;w)=\codim(X_w)-\codim(X_u)-\codim(X_v)$,
one easily shows that $c_{u,v}^w=a_{u,v}^w$ if $N(u,v;w)=0$, whereas 
$c_{u,v}^w=0$ if $N(u,v;w)<0$.

\smallskip

In the case where $X$ is a Grassmannian, Buch obtained a
combinatorial description of the structure constants $c_{u,v}^w$
which implies that they have alternating signs:
$$
(-1)^{N(u,v;w)}\, c_{u,v}^w\geq 0
$$
for all $u$, $v$, $w$. And he conjectured that the latter result holds
for all flag varieties (\cite{Bu} Conjecture 9.2).

\smallskip

This conjecture is proved in the present paper. In fact it is a direct
consequence of the following result, which answers a question of
W.~Graham.

\begin{theorem}\label{signs}
Let $Y$ be a closed subvariety of a complex flag variety $X$; write
$$
[\cO_Y]=\sum_w\, c_Y^w\,[\cO_{X_w}]
$$
in $K(X)$. If $Y$ has rational singularities (e.g., if $Y$ is
nonsingular), then the coefficients $c_Y^w$ satisfy 
$$
(-1)^{\codim(X_w)-\codim(Y)} \, c_Y^w\geq 0.
$$ 
\end{theorem}

Using the duality involution of $K(X)$, this may be reformulated in a
more appealing way: {\it the class of the dualizing sheaf $\omega_Y$
of any subvariety $Y$ having rational singularities, is a nonnegative 
combination of classes of dualizing sheaves of Schubert varieties}
(the latter classes form another natural basis of $K(X)$).

\smallskip

To deduce Buch's conjecture from Theorem \ref{signs}, one shows that
$[\cO_{X_u}]\cdot[\cO_{X_v}]=[\cO_Y]$, where
$Y$ is the intersection of $X_u$ with a general translate
of $X_v$ ; furthermore, since Schubert varieties have rational
singularities, the same holds for $Y$ (see Lemma \ref{general}).

\smallskip

We now sketch a proof of Theorem \ref{signs} in the simplest case,
where $X$ is a projective space and $Y$ is nonsingular. The Schubert
varieties in $X=\mP^n$ are a flag of linear subspaces $\mP^m$, where
$0\leq m \leq n$. Writing 
$$
[\cO_Y]=\sum_{m=0}^n\, c_Y^m \,[\cO_{\mP^m}]
$$
and noticing that the Euler characteristic $\chi(\cO_{\mP^m}(-1))$
vanishes for $m>0$, we see that $c_Y^0=\chi(\cO_Y(-1))$. By the
Kodaira vanishing theorem, $H^i(\cO_Y(-1))=0$ for $i<\dim(Y)$, so that
$$
(-1)^{\dim(Y)} \,c_Y^0\geq 0.
$$ 
More generally, one obtains 
$$
c_Y^m =\chi(\cO_{Y\cap \mP^{n-m}}(-1))
$$ 
for any linear subspace $\mP^{n-m}$ in general position with respect
to $Y$; then $Y\cap \mP^{n-m}$ is nonsingular, and the same argument
yields 
$$
(-1)^{\dim(Y)-m} \, c_Y^m\geq 0.
$$

\smallskip

This argument adapts to subvarieties $Y$ having rational
singularities, by Grauert-Riemenschneider's generalization of the
Kodaira vanishing theorem. It also shows that the statement of Theorem
\ref{signs} does not hold for all closed subvarieties of $\mP^n$. For
instance, one may check that the projective cone $Y$ over a
nonsingular rational curve of degree $d$ in $\mP^{n-1}$ satisfies
$c_Y^0 \leq n-d$. This yields examples of (singular) projective
surfaces $Y$ with arbitrarily negative $c_Y^0$. 

\smallskip

However, Theorem \ref{signs} may be extended to all closed
subvarieties of complex flag varieties, by replacing $[\cO_Y]$
with $\sum_{i=0}^{\dim(Y)}\,(-1)^i\,[R^i\varphi_*\cO_{\tY}]$ where
$\varphi:\tY\to Y$ is a desingularization (see the Remark at the end
of Section 2). On the other hand, we do not know whether Theorem
\ref{signs} extends to positive characteristics, already for
nonsingular surfaces in projective space: although the Kodaira
vanishing theorem does not hold in this setting, $\chi(\cO_Y(-1))$ may
well be nonnegative (see \cite{Ray}).

\smallskip

This article is organized as follows. Section 1 gathers preliminary
results concerning products of classes of structure sheaves and
dualizing sheaves in the Grothendieck ring of flag varieties. 
In Section 2, we generalize the decomposition
$$
[\cO_Y]=\sum_{m=0}^n \,\chi(\cO_{Y\cap \mP^{n-m}}(-1)) \,[\cO_{\mP^m}]
$$
to any Cohen-Macaulay subvariety $Y$ of a flag variety $X$. For this,
we construct a degeneration in $X\times X$ of the diagonal of $Y$,
which is interesting in its own right (Theorem \ref{degeneration}). 
In Section 3, we obtain an analogue of the vanishing theorem
$$
H^i(\cO_{Y\cap\mP^{n-m}}(-1))=0 \text{ for } i<\dim(Y)-m,
$$
where $Y$ has rational singularities ; this completes the proof of
Theorem \ref{signs}. 

\smallskip

In the final Section 4, we adapt these arguments to prove another
positivity result: the class of the restriction of any globally
generated line bundle to any Schubert variety, has nonnegative
coefficients in the basis of classes of Schubert structure
sheaves. That result was first obtained by W.~Fulton and A.~Lascoux
for the variety of complete flags, using the combinatorics of
Grothendieck polynomials (see \cite{FL}). It was generalized to all
flag varieties by H.~Pittie and A.~Ram (see \cite{PR}) using the
Littelmann path method, and also by O.~Mathieu (see \cite{Mat}) using
representation theory. Our geometric approach expresses the
coefficients as the dimensions of certain spaces of sections of the
line bundle (Theorem \ref{line}). 

\smallskip

\noindent
{\it Acknowledgements}. Many thanks to W.~Graham, S.~Kumar and
E.~Peyre for useful discussions and suggestions.

\medskip

\noindent
{\bf Notation and conventions.} The ground field is the field
$\mC$ of complex numbers. An equidimensional reduced scheme of
finite type over $\mC$ will be called a {\it variety} ; with this
convention, varieties need not be irreducible.

A {\it desingularization} of a variety $Y$ is a nonsingular
variety $\tY$ together with a proper birational morphism 
$\varphi:\tY\to Y$. The singularities of $Y$ are {\sl rational} 
if $Y$ is normal and if there exists a desingularization
$\varphi:\tY\to Y$ such that $R^i\varphi_*(\cO_{\tY})=0$ for
all $i\geq 1$. Equivalently, $Y$ is Cohen-Macaulay and the natural map
$\varphi_*\omega_{\tY}\to\omega_Y$ is an isomorphism, where
$\omega_{\tY}$  (resp. $\omega_Y$) denotes the dualizing sheaf of
$\tY$ (resp. $Y$). If $Y$ has rational singularities, then 
the preceding conditions hold for any desingularization 
(for these results, see e.g. \cite{KKMS} p.~50). Furthermore, since
$Y$ is normal, $\omega_Y$ is the double dual of the sheaf
$\wedge^{\dim(Y)}\Omega_Y^1$, that is, the direct image of the
sheaf of differential forms of top degree on the nonsingular locus
$Y^{\reg}$.

We next turn to notation concerning flag varieties.
Let $G$ be a simply connected semisimple algebraic group. Choose
opposite Borel subgroups $B$ and $B^-$, with common torus $T$ ; let
$\cX(T)$ be the group of characters of $T$, also called weights. In
the root system $\Phi$ of $(G,T)$, we have the subset $\Phi^+$ of
positive roots (that is, of roots of $(B,T)$), and the subset
$\Delta=\{\alpha_1,\ldots,\alpha_r\}$ of simple roots. Let
$\omega_1,\ldots,\omega_r$ be the corresponding fundamental weights ;
they form a basis of $\cX(T)$. Let $\rho=\omega_1+\cdots+\omega_r$,
this equals the half sum of positive roots. 

We also have the Weyl group $W$ of $(G,T)$, generated by the simple
reflections $s_1,\ldots,s_r$ corresponding to the simple roots. This
defines the length function $\ell$ and the Bruhat order $\leq$ on
$W$. Let $w_o$ be the longest element of $W$, then $B^-=w_o B w_o$.

Let $P$ be a parabolic subgroup of $G$ containing $B$ and let $W_P$
be the Weyl group of $(P,T)$, a parabolic subgroup of $W$ ; let
$w_{o,P}$ be the longest element of $W_P$.
Each right $W_P$-coset in $W$ contains a unique element of minimal
length; this defines the subset $W^P$ of minimal representatives of
the quotient $W/W_P$. This subset is invariant under the map 
$w\mapsto w_oww_{o,P}$ ; the induced bijection of $W^P$ reverses
the Bruhat order. Notice that $W^P=W$ if and only if $P=B$.

The homogeneous variety $X=G/P$ is called a flag variety, the
full flag variety being $G/B$. For any weight $\lambda$ regarded as 
a character of $B$, let $\cL_{G/B}(\lambda)$ be the corresponding
$G$-linearized invertible sheaf on $G/B$. The assignment
$\lambda\mapsto\cL_{G/B}(\lambda)$ yields an isomorphism from $\cX(T)$
to the Picard group of $G/B$; the dominant weights correspond to the
globally generated invertible sheaves. The dualizing
sheaf of $G/B$ is $\cL_{G/B}(-2\rho)$. Via pullback under the natural
map $G/B\to G/P$, the Picard group of $G/P$ identifies with the
subgroup of $\cX(T)$ consisting of restrictions of characters of $P$.

The $T$-fixed points in $G/P$ are the $e_{wP}=wP/P$ ($w\in W/W_P$); 
we index them by $W^P$. The $B$-orbit $C_{wP}=B e_{wP}$ is a Bruhat
cell, isomorphic to affine space of dimension $\ell(w)$; its closure
in $X$ is the Schubert variety $X_{wP}$. The complement 
$X_{wP} - C_{wP}$ is the boundary $\partial X_{wP}$ ; it has pure
codimension $1$ in $X_{wP}$. 

We shall also need the opposite Bruhat cell $C^-_{wP}=B^- e_{wP}$ of
codimension $\ell(w)$ in $X$, the opposite Schubert variety
$X^{wP}=\overline{C^-_{wP}}$, and its boundary $\partial X^{wP}$. 
Then 
$$
X^{wP}=w_o X_{w_ow w_{o,P}P}.
$$ 
For any $v$ and $w$ in $W^P$, we have 
$v\leq w\Leftrightarrow X_{vP}\subseteq X_{wP} \Leftrightarrow
X^{wP}\subseteq X^{vP}$. 
Equivalently,
$$
\partial X_{wP}=\bigcup_{v\in W^P,\,v< w} X_{vP} ~\text{ and }~
\partial X^{wP}=\bigcup_{v\in W^P,\,v> w} X^{vP}.
$$

By \cite{Ram} Proposition 2 and Theorem 4, all Schubert varieties have
rational singularities; in particular, they are normal and
Cohen-Macaulay. If in addition $X$ is the full flag variety, then we 
use the simpler notation $e_w$, $C_w$, $X_w$ ,$\ldots$. We then have
$$
\omega_{X_w}=\cL_{X_w}(-\rho)(-\partial X_w).
$$

\section{Preliminaries on Grothendieck groups}

For an arbitrary nonsingular variety $X$, let $K(X)$ be the Grothendieck
group of the category of coherent sheaves on $X$ ; the class in $K(X)$
of a coherent sheaf $\cF$ will be denoted $[\cF]$. Recall that $K(X)$
is isomorphic to the Grothendieck group of vector bundles over $X$. 
The tensor product of vector bundles defines a product $\cdot$ on $K(X)$
making it a commutative ring with unit ; the duality of vector bundles
defines an involutive automorphism $*$ of that ring. For any coherent
sheaves $\cF$, $\cG$ on $X$, we have
$$
[\cF] \cdot [\cG] = \sum_{i = 0}^{\dim(X)} \, (-1)^i\,
[Tor_i^X(\cF,\cG)]
\text{ and }
[\cF]^* =  \sum_{i = 0}^{\dim(X)} \, (-1)^i\,
[Ext^i_X(\cF,\cO_X)].
$$
If in addition $X$ is complete, then the Euler characteristic
$$
\cF\mapsto\chi(\cF)=\sum_{i=0}^{\dim(X)}\, (-1)^i\,h^i(\cF)
$$ 
yields the pushforward map $\chi:K(X)\to\mZ$.

We associate to any closed subscheme $Y$ of $X$ the class $[\cO_Y]$ 
of its structure sheaf. If in addition $Y$ is Cohen-Macaulay and
equidimensional, its dualizing sheaf
$$
\omega_Y = Ext_X^{\codim(Y)}(\cO_Y,\omega_X)
$$
yields another class $[\omega_Y]$, and we have
$$
[\cO_Y]^* = (-1)^{\codim(Y)} [\omega_Y]\cdot[\omega_X]^*.
$$

The product of two such classes associated with subschemes ``in general
position'' is determined by the following result, a variant of a lemma
of Fulton and Pragacz (\cite{FP}, p.~108).

\begin{lemma}\label{intersection}
Let $Y$ and $Z$ be closed subschemes of a nonsingular variety $X$. If
$Y$ and $Z$ are equidimensional, Cohen-Macaulay and intersect properly
in $X$, then their scheme-theoretic intersection $Y\cap Z$ is
equidimensional and Cohen-Macaulay as well. Furthermore,
$$
\omega_{Y\cap Z}=
\omega_Y\otimes_{\cO_X}\omega_Z\otimes_{\cO_X}\omega_X^{-1},
$$ 
and $Tor_i^X(\cO_Y,\cO_Z) = 0 = Tor_i^X(\omega_Y,\omega_Z)$ 
for all $i\geq 1$. As a consequence,
$$
[\cO_Y]\cdot[\cO_Z]=[\cO_{Y\cap Z}] \text{ and }
[\omega_Y]\cdot[\omega_Z]=[\omega_{Y\cap Z}]\cdot[\omega_X].
$$
\end{lemma}

\begin{proof}
Notice that $Y\times Z$ intersects properly the diagonal, $\diag(X)$
in $X\times X$, along $\diag(Y\cap Z)$. Thus, a sequence of local
equations of $\diag(X)$ in $X\times X$ at any point of 
$\diag(Y\cap Z)$ restrict to a regular sequence in the local ring of 
$Y\times Z$ at that point. Since $Y\times Z$ is Cohen-Macaulay, it
follows that the same holds for $\diag(Y\cap Z)$, and that
$Tor_i^{X\times X}(\cO_{Y\times Z},\cO_{\diag(X)})=0$ 
for any $i\geq 1$.

Let $c$ (resp. $d$) denote the codimension of $Y$ (resp. $Z$) in
$X$. Then we may choose a locally free resolution $\cL$ (resp. $\cM$)
of the sheaf of $\cO_X$-modules $\cO_Y$ (resp. $\cO_Z)$ of length $c$
(resp. $d$). Now $\cL\otimes_{\mC}\cM$ is a locally free resolution of
the sheaf of $\cO_{X\times X}$-modules $\cO_{Y\times Z}$. By the
preceding step, it follows that 
$$
\cL\otimes_{\cO_X}\cM \cong 
(\cL\otimes_{\mC}\cM)\otimes_{\cO_{X\times X}}\cO_{\diag(X)}
$$ 
is a locally free resolution of $\cO_{Y\cap Z}$. Thus,
$Tor_i^X(\cO_Y,\cO_Z)=0$ for all $i\geq 1$.

Since $Y$ is Cohen-Macaulay, the dual complex $\cL^*$ is a locally
free resolution of the sheaf
$$
Ext_X^c(\cO_Y,\cO_X)=\omega_Y\otimes_{\cO_X}\omega_X^{-1}.
$$
Likewise, $\cM^*$ (resp. $(\cL\otimes\cM)^*=\cL^*\otimes\cM^*$) is a
locally free resolution of $\omega_Z\otimes_{\cO_X}\omega_X^{-1}$
(resp. of $\omega_{Y\cap Z}\otimes_{\cO_X}\omega_X^{-1}$). It follows
that 
$$
\omega_Y\otimes_{\cO_X} \omega_X^{-1}\otimes_{\cO_X} 
\omega_Z\otimes_{\cO_X}\omega_X^{-1}
\cong \omega_{Y\cap Z}\otimes_{\cO_X} \omega_X^{-1}
$$
and that $Tor_i^X(\omega_Y,\omega_Z)=0$ for all $i\geq 1$.
\end{proof}

We next turn to the case where $X$ is a flag variety ; we shall obtain
variants of Kleiman's transversality theorem, see \cite{Kl}. In
what follows, a statement holds ``for general $g\in G$'' if it holds
for all $g$ in a nonempty open subset of $G$.

\begin{lemma}\label{general}
Let $Y$ be a closed subvariety of $X=G/P$ and let $w\in W^P$. Then the
translate $g X_{wP}$ intersects properly $Y$, for general $g\in G$. 

If in addition $Y$ is Cohen-Macaulay, then $Y\cap g X_{wP}$ is
Cohen-Macaulay as well, for general $g$. As a consequence, 
$$
[\cO_Y]\cdot[\cO_{X_{wP}}]=[\cO_{Y\cap g\cdot X_{wP}}] \text{ and }
[\omega_Y]\cdot[\omega_{X_{wP}}]=[\omega_{Y\cap g\cdot X_{wP}}]\cdot[\omega_X].
$$
If in addition $Y$ is normal, then $Y\cap gX_{wP}$ is normal as well.

Finally, if $Y$ has rational singularities, then $Y\cap g X_{wP}$
has rational singularities as well. 
\end{lemma}

\begin{proof}
The first assertion follows from Kleiman's transversality theorem ; we
recall the proof, since we shall repeatedly use its ingredients. Let
$$
i:Y\to X
$$ 
be the inclusion and let 
$$
m:G\times X_{wP}\to X,~(g,x)\mapsto g x
$$ 
be the ``multiplication'' map. Notice that $G$ acts on $G\times X_{wP}$
via left multiplication on $G$, and that $m$ is a $G$-equivariant
morphism to $G/P$. Thus, $m$ is a locally trivial fibration with fiber 
$m^{-1}(P/P)$. The latter is isomorphic to
$\overline{Pw^{-1}B}$, the pullback in $G$ of a Schubert variety.
Therefore, the fiber of $m$ has rational singularities. 

Now consider the cartesian product
$$
Z=Y\times_X (G\times X_{wP})
$$ 
with projections $\iota$ to $G\times X_{wP}$, and $\mu$ to $Y$. Let
$p:G\times X_{wP}\to G$ be the projection and let
$\pi=p\circ\iota$ as displayed in the following commutative diagram:
 
$$
\CD
 G @<{\pi}<< Z @>{\mu}>> Y\\
 @V{id}VV    @V{\iota}VV @V{i}VV\\
 G @<{p}<<   G\times X_{wP} @>{m}>> X
\endCD
$$
By definition, the square on the right is cartesian, so that $\mu$ is
also a locally trivial fibration with fiber $\overline{Pw^{-1}B}$.  
As a consequence, $Z$ is a variety, and we have  
$$
\dim(Z)=\dim(Y)+\dim(G\times X_{wP})-\dim(X).
$$ 
Furthermore, the morphism $\pi$ is proper, and its fiber at each $g\in G$
identifies with $Y\cap g X_{wP}$. Thus, for general $g$, the intersection
$Y\cap g X_{wP}$ is either empty or equidimensional of dimension
$\dim(Y)+\dim(X_{wP})-\dim(X)$. This proves the first assertion.

If in addition $Y$ is Cohen-Macaulay, then by Lemma \ref{intersection},
the intersection $Y\cap g X_{wP}$ is Cohen-Macaulay whenever it is
proper. If in addition $Y$ is normal, then $Y\cap gX_{wP}$ is nonsingular
in codimension $1$, for general $g$ (since 
$Y^{\reg}\cap gX_{wP}^{\reg}$ is nonsingular). Therefore, 
$Y\cap gX_{wP}$ is normal, by Serre's criterion. 

Assume now that $Y$ has rational singularities ; then $Z$ has rational
singularities as well. Now the following easy result completes the
proof of Lemma \ref{general}.

\begin{lemma}\label{fibres}
Let $Z$ and $S$ be varieties and let $\pi:Z\to S$ be a morphism. If
$Z$ has rational singularities, then the same holds for the general
fibers of $\pi$.
\end{lemma}

\begin{proof}
We may reduce to the situation where $Z$ is affine, $S$ is nonsingular
and $\pi$ is flat with connected fibers. Let $F$ be a general fiber of
$\pi$ ; then $F$ is Cohen-Macaulay, since $Z$ is.

Choose a desingularization $\varphi:\tZ\to Z$ and let
$\tilde\pi=\pi\circ\varphi$. Then $\tilde F=\varphi^{-1}(F)$ is a
general fiber of $\tilde\pi$, and hence is connected. By generic
smoothness, $\tilde F$ is nonsingular, so that
$\varphi$ restricts to a desingularization 
$\psi:\tilde F\to F$. Since $Z$ has rational singularities, the
map $\varphi_*\omega_{\tZ}\to \omega_Z$ is an isomorphism ; since 
$\omega_F$ is the restriction of $\omega_Z\otimes\pi^*\omega_S^{-1}$
to $F$, and similarly for $\omega_{\tilde F}$,
it follows that the map $\psi_*\omega_{\tilde F}\to \omega_F$ is an
isomorphism as well.
\end{proof}
\end{proof}

Consider, for example, $Y=X^{vP}$ where $v\in W^P$. Since $Y$
(resp. $X_{wP}$) is invariant under $B^-$ (resp. $B$) and the product
$B^- B$ is open in $G$, we see that $X^{vP}\cap X_{wP}$ is the
intersection of $X^{vP}$ with a general translate of $X_{wP}$, and
hence has rational singularities. 

In fact $X^{vP}\cap X_{wP}$ is the closure of $C^{vP}\cap C_{wP}$, and
hence is irreducible of dimension $\ell(w)-\ell(v)$ (or empty), see
\cite{Ri} Theorem 3.7. Considering $T$-fixed points, one sees that
$X^{vP}\cap X_{wP}$ is nonempty if and only if $v\leq w$. 

If in addition $X$ is the full flag variety, then we obtain
$$
\omega_{X^v\cap X_w}=\cO_{X^v\cap X_w}
(-(X_w\cap \partial X^v) -(X^v\cap\partial X_w)).
$$
Likewise, $X_{vP} \cap w_oX_{wP} = X_{vP}\cap X^{w_oww_{o,P}P}$ is the
intersection of $X_{vP}$ with a general translate of $X_{wP}$.

\medskip

As a final preliminary result, we study the boundaries of Schubert
varieties:

\begin{lemma}\label{boundary}
For every $w\in W^P$, the boundary $\partial X_{wP}$ is
Cohen-Macaulay. Moroever, for every Cohen-Macaulay subvariety $Y$ of
$X$, the intersection $Y\cap gX_{wP}$ is Cohen-Macaulay for general
$g\in G$. As a consequence, 
$$
Ext^i_{Y\cap gX_{wP}}(\cO_{Y\cap gX_{wP}}
(-Y\cap g\partial X_{wP}),\omega_{Y\cap gX_{wP}})=0 
\text{ for any }i\geq 1.
$$ 
If in addition $Y$ is normal, then 
$$
Hom(\cO_{Y\cap gX_{wP}}(-Y\cap g\partial X_{wP}),\omega_{Y\cap gX_{wP}})
=\omega_{Y\cap gX_{wP}}(Y\cap g\partial X_{wP}).
$$
\end{lemma}

\begin{proof}
If $X=G/B$ then $\cO_{X_w}(-\partial X_w)$ is locally isomorphic to
$\omega_{X_w}$, so that the ideal sheaf of $\partial X_w$ in $X_w$ is
Cohen-Macaulay. Since $X_w$ is also Cohen-Macaulay and $\partial X_w$
has pure codimension $1$, it follows easily that $\partial X_w$ is
Cohen-Macaulay as well.

In the general case where $X=G/P$, notice that the natural map 
$G/B\to G/P$ restricts to a proper surjective morphism 
$$
\eta:X_w\to X_{wP}
$$
which maps isomorphically $C_w$ to $C_{wP}$ (since $w\in W^P$). Thus,
$\eta^{-1}(\partial X_{wP})=\partial X_w$ (as sets). By \cite{Ram}
Theorem 2 and Proposition 3, we have 
$$
\eta_*\cO_{X_w}=\cO_{X_{wP}}~\text{ and }~
R^i\eta_*(\cO_{X_w})=0 ~\text{ for }~i\geq 1.
$$ 
It follows that 
$\eta_*\cO_{X_w}(-\partial X_w)=\cO_{X_{wP}}(-\partial X_{wP})$.
We claim that
$$
R^i\eta_*(\cO_{X_w}(-\partial X_w))=0 \text{ for } i\geq 1.
$$

To see this, choose a reduced decomposition of $w$ and let 
$$
\psi:\tX_w\to X_w
$$ 
be the corresponding standard desingularization (see \cite{Ram}). Then
$B$ acts in $\tX_w$, and $\psi$ is $B$-equivariant. Furthermore, $\tX_w$
contains a dense $B$-orbit, mapped isomorphically to $C_w$ by $\psi$ ;
the complement $\partial\tX_w$ of this orbit is a union of nonsingular 
irreducible divisors $\tX_1,\ldots,\tX_{\ell}$ intersecting 
transversally. By \cite{Ram} Proposition 2, we have 
$$
\omega_{\tX_w}=
\cO_{\tX_w}(-\partial\tX_w)\otimes\psi^*\cL_{X_w}(-\rho).
$$
Let 
$$
\tilde \eta =  \eta\circ \psi: \tX_w\to X_w,
$$
this is a desingularization of $X_w$. Since $\cL_{X_w}(\rho)$ is
ample, we have 
$$
R^i\tilde\eta_*(\omega_{\tilde X_w}\otimes\psi^*\cL_{X_w}(\rho))=0
$$
for $i\geq 1$, by the Grauert-Riemenschneider vanishing theorem
(see \cite{GR}) ; furthermore,
$$
R^i\psi_*(\omega_{\tilde X_w}\otimes\psi^*\cL_{X_w}(\rho))=
R^i\psi_*(\omega_{\tilde X_w})\otimes\cL_{X_w}(\rho)
$$
vanishes for $i\geq 1$ as well, and 
$$
\psi_*(\omega_{\tilde X_w}\otimes\psi^*\cL_{X_w}(\rho))=
\psi_*\omega_{\tilde X_w}\otimes\cL_{X_w}(\rho)=
\omega_{X_w}\otimes\cL_{X_w}(\rho)=
\cO_{X_w}(-\partial X_w).
$$
Now the Leray spectral sequence for $\tilde\eta=\eta\circ\psi$ implies
the claimed vanishing.

Using that claim and duality for the morphism $X_w\to G/P$, we obtain
$$\displaylines{
Ext^i_X(\cO_{X_{wP}}(-\partial X_{wP}),\omega_X)
=R^{i-\codim(X_{wP})}\eta_*(\omega_{X_w}(\partial X_w))
\hfill\cr\hfill
=R^{i-\codim(X_{wP})}\eta_*(\cL_{X_w}(-\rho)).
\cr}$$
Thus, to prove that the sheaf $\cO_{X_{wP}}(-\partial X_{wP})$ is
Cohen-Macaulay, it suffices to check the vanishing of
$R^i\eta_*(\cL_{X_w}(-\rho))$ for $i\geq 1$. We deduce this from the
Kawamata-Viehweg vanishing theorem (see e.g. \cite{EV}) as follows.

Notice that $\partial \tX_w$ is the support of a very ample divisor of
$\tX_w$ (to see this, consider a $B$-linearized very ample invertible
sheaf $\cM$ on $\tX_w$ and a $B$-semi-invariant section $\sigma$ of
$\cM$ that vanishes on $\partial\tX_w$ ; then the zero set of $\sigma$
is exactly $\partial \tX_w$). 
Thus, we may choose positive integers $b_1,\ldots,b_{\ell}$ such that
the divisor $b_1\tX_1+\cdots+b_{\ell}\tX_{\ell}$ is very
ample. Choose also a positive integer $N>\max(b_1,\ldots,b_{\ell})$
and let $a_1=N-b_1,\ldots,a_{\ell}=N-b_{\ell}$. Finally, let
$\cL=\cO_{\tX_w}(\tX_w)$ and $D=a_1\tX_1+\cdots+a_{\ell}\tX_{\ell}$. 
Then the invertible sheaf
$$
\cL^N(-D)=\cO_{\tX_w}(b_1\tX_1+\cdots+b_{\ell}\tX_{\ell})
$$
is very ample. By \cite{EV} Corollary 6.11, it follows that 
$R^i\tilde \eta_*(\omega_{\tX_w}(\tX_w))=0$ for $i\geq 1$, that is,
$R^i\tilde \eta_*(\psi^*\cL_{X_w}(-\rho))=0$. As above, the Leray
spectral sequence for $\tilde\eta=\eta\circ \psi$ yields the vanishing
of $R^i\eta_*(\cL_{X_w}(-\rho))$ for $i\geq 1$ ; thus,
$\cO_{X_{wP}}(-\partial X_{wP})$, and hence $\partial X_{wP}$, is
Cohen-Macaulay.
  
By Lemma \ref{general}, the same holds for $Y\cap g\partial X_{wP}$ for
general $g$. Since $Y\cap g\partial X_{wP}$ has pure codimension $1$ in
$Y\cap gX_{wP}$, and the latter is Cohen-Macaulay, the ideal sheaf
$\cO_{Y\cap g X_{wP}}(-Y\cap g\partial X_{wP})$ is Cohen-Macaulay as well.
This implies the vanishing of 
$Ext^i_{Y\cap gX_{wP}}(\cO_{Y\cap gX_{wP}}
(-Y\cap g\partial X_{wP},\omega_{Y\cap gX_{wP}})$ for $i\geq 1$.
If in addition $Y$ is normal, then so is $Y\cap gX_{wP}$ by Lemma
\ref{general}. Thus, the sheaf $\omega_{Y\cap gX_{wP}}$ is reflexive of
rank $1$ ; this implies the latter assertion.
\end{proof}

\section{A degeneration of the diagonal}

Let $Y$ be a closed subvariety of a flag variety $X$. Consider the
diagonal, $\diag(Y)\subseteq X\times X$. We shall construct a flat
degeneration of $\diag(Y)$ to a union of products of subvarieties, by
taking limits under all ``positive'' one-parameter subgroups. 

Specifically, let $T$ act linearly in affine space $\mA^r$ with
weights $-\alpha_1,\ldots,-\alpha_r$. Then the $T$-orbit of the point
$(1,\ldots,1)$ is the complement of the union of all coordinate
hyperplanes, and the isotropy group of that point is the center
$Z(G)$. Thus, the orbit $T\cdot(1,\ldots,1)$ is isomorphic to the
adjoint torus $T_{\ad}=T/Z(G)$.

Let $\cY$ be the closure in $X\times X\times\mA^r$ of the subset
$$
\{(ty,y,\alpha_1(t),\ldots,\alpha_r(t))~\vert~ y\in Y, ~t\in T\}
$$
and let
$$
\pi_Y:\cY\to \mA^r,~p_Y:\cY\to X\times X
$$
be the projections. Clearly, the morphism $\pi_Y$ is proper and its
fibers identify via $p_Y$ with closed subschemes of $X\times X$. This
yields an isomorphism
$$
\pi_Y^{-1}(1,\ldots,1)\cong \diag(Y).
$$

Notice that $T$ acts in $X\times X\times \mA^r$ by
$$
t\cdot(x_1,x_2,t_1,\ldots,t_r)=
(tx_1,x_2,\alpha_1(t)t_1,\ldots,\alpha_r(t)t_r)
$$
and leaves $\cY$ invariant ; furthermore, $\pi_Y$ is equivariant.
Thus, $\pi_Y$ is surjective, and restricts to a trivial fibration over
the orbit $T\cdot(1,\ldots,1)$ with fiber $\diag(Y)$. Some less
obvious properties of $\pi_Y$ are summarized in the following
statement.

\begin{theorem}\label{degeneration}
Let $Y$ be a Cohen-Macaulay subvariety of $X$, such that the
intersection $Y\cap X^{wP}$ is proper and reduced for all $w\in W^P$.
Then, with preceding notation, $\cY$ is Cohen-Macaulay ; furthermore,
$\pi_Y$ is flat with reduced fibers, and
$$
\pi_Y^{-1}(0,\ldots,0)=\bigcup_{w\in W^P} X_{wP}\times (Y\cap X^{wP}).
$$
\end{theorem}

\begin{proof}
We begin with the case where $Y=X$ ; we then set $\cY=\cX$. We show
how to obtain $\pi_X:\cX\to\mA^r$ by base change from a degeneration of
$\diag(X)$ constructed in \cite{BP} \S 7. 

Let $G_{\ad}=G/Z(G)$ be the adjoint group of $G$ and let
$\oG_{\ad}$ be its wonderful completion ; this is a
nonsingular projective variety where $G\times G$ acts with
a dense orbit isomorphic to 
$(G\times G)/(Z(G)\times Z(G))\diag(G)\cong G_{\ad}$, and a
unique closed orbit isomorphic to $(G\times G)/(B\times B^-)$. Let
$P_{\ad}=P/Z(G)$ and let $\oP$ be its closure in
$\oG_{\ad}$. Since $\oP_{\ad}$ is invariant under
the action of $P\times P$, we may form the associated fiber bundle
$$
p:G\times G\times^{P\times P}\oP_{\ad}\to G/P\times G/P
=X\times X.
$$
On the other hand, the map 
$G\times G\times\oP_{\ad}\to\oG_{\ad}$, 
$(x,y,z)\mapsto (x,y) z$ factors through a map
$$
\pi:G\times G\times^{P\times P}\oP_{\ad}\to 
\oG_{\ad},
$$
which is clearly surjective and $G\times G$-equivariant. Furthermore,
the product map 
$$
p\times \pi: G\times G\times^{P\times P}\oP_{\ad}\to 
X\times X\times\oG_{\ad}
$$
is a closed immersion, with image the ``incidence variety''
$$
\{(xP,yP,z)~\vert~ (x,y)\in G\times G, z\in (x,y)\oP_{\ad}\}.
$$

By \cite{BP} \S 7, $\pi$ is flat, with reduced Cohen-Macaulay fibers ;
these identify with closed subschemes of $X\times X$ via $p$. The
fiber at the identity element of $G_{\ad}$ (resp. at the unique 
$(B\times B^-)$-fixed point $z$) identifies with $\diag(X)$
(resp. $\cup_{w\in W^P}\,X_{wP}\times X^{wP}$). Furthermore, the
closure $\oT_{\ad}$ of the torus $T_{\ad}$ in $\oG_{\ad}$ is
a nonsingular $T\times T$-equivariant completion of that torus,
containing $z$ as a fixed point. 

Let $\oT_{\ad,}z$ be the unique $T\times T$-invariant open
affine neighborhood of $z$ in $\oT_{\ad}$. Then
$\oT_{\ad,z}$ is equivariantly isomorphic to affine $r$-space
where $T\times T$ acts linearly with weights 
$(-\alpha_1,\alpha_1),\ldots,(-\alpha_r,\alpha_r)$. Thus, for the
action of $T$ by left multiplication, $\oT_{\ad,z}$ is
isomorphic to $\mA^r$.

We claim that the subvariety $\pi^{-1}(\oT_{\ad,z})$ of 
$X\times X\times \mA^r$ equals $\cX$. To see this, note that
$$\displaylines{
\pi^{-1}(T)=\{(xP,yP,z)~\vert~ (x,y)\in G\times G, 
z\in T_{\ad}\cap(x,y)\oP_{\ad}\}
\hfill\cr\hfill
=\{(zyP,yP,z)~\vert~ y\in G,z\in T_{\ad}\}
=\{(t\xi,\xi,t\cdot(1,\ldots,1))~\vert~ \xi\in X, ~t\in T\}
\cr}$$
since 
$T_{\ad}\cap (x,y)\oP_{\ad}=T_{\ad}\cap xP_{\ad}y^{-1}$. 
It follows that $\pi^{-1}(\oT_{\ad,z})$ contains $\cX$ as an
irreducible component. Furthermore, $\pi$ restricts to a flat morphism 
from the complement $\pi^{-1}(\oT_{\ad,z}) - \cX$, to
$\oT_{\ad,z}$. If this complement is not empty, then its
image meets the open subset $T_{\ad}$ of $\oT_{\ad,z}$, a
contradiction. This proves the claim, and hence all assertions of
Theorem \ref{degeneration} in the case where $Y=X$.

In the general case, we consider
$$
\cY'=\cX\cap (X\times Y\times \mA^r)
$$
(scheme-theoretical intersection in $X\times X\times \mA^r$), with
projection $\pi'_Y:\cY'\to\mA^r$. Notice that $\cY$ is contained in
$\cY'$, and that 
$$
\cY\cap\pi_Y^{-1}(T\cdot(1,\ldots,1))=
\cY'\cap\pi^{'-1}_Y(T\cdot(1,\ldots,1)).
$$
Thus, $\cY$ is an irreducible component of $\cY'$ ; the latter is
invariant under the action of $T$ in $X\times X\times \mA^r$. 

We claim that $\cX$ intersects properly $X\times Y\times \mA^r$ in
$X\times X\times \mA^r$, that is, every irreducible component $C$ of
$\cY'$ has dimension equal to
$$
\dim(\cX)+\dim(X\times Y\times\mA^r)-\dim(X\times X\times \mA^r)
= \dim(Y)+r.
$$
In fact, it suffices to check that $\dim(C)\leq \dim(Y)+r$. Since $C$
is $T$-invariant and $T$ acts attractively in $\mA^r$ with fixed point 
$(0,\ldots,0)$, it suffices in turn to show that 
$\dim(C\cap\pi^{'-1}_Y(0,\ldots,0))\leq \dim(Y)$. But
$$\displaylines{
C\cap\pi^{'-1}_Y(0,\ldots,0)\subseteq 
\cX\cap(X\times Y\times \{(0,\ldots,0)\}
\hfill\cr\hfill
\cong(\bigcup_{w\in W^P} X_{wP}\times X^{wP})\cap (X\times Y)=
\bigcup_{w\in W^P} X_{wP}\times(Y\cap X^{wP}).
\cr}$$
And the latter is equidimensional of dimension $\dim(Y)$, since
$Y$ intersects properly all $X^{wP}$. This proves our claim.

Since $\cX$ and $X\times Y\times\mA^r$ are Cohen-Macaulay subvarieties
of $X\times X\times\mA^r$ intersecting properly, then $\cY'$ is
equidimensional and Cohen-Macaulay, by Lemma
\ref{intersection}. Furthermore, the morphism $\pi'_Y:\cY'\to \mA^r$ is
equidimensional by the proof of the preceding claim ; therefore,
$\pi'_Y$ is flat. As in the first step of the proof, it follows that
$\cY'$ equals $\cY$ as sets. Furthermore,
$$
\pi^{'-1}_Y(T\cdot(1,\ldots,1))=
\cX\cap (X\times Y\times T\cdot(1,\ldots,1))
$$
is clearly reduced, so that $\cY'$ is generically reduced. Since it is
Cohen-Macaulay, it is reduced, and $\cY'$ equals $\cY$ as subschemes.

Likewise, the fiber $\pi^{-1}_Y(0,\ldots,0)$ equals 
$\cup_{w\in W^P} X_{wP}\times (Y\cap X^{wP})$ as sets. Furthermore, this
fiber is generically reduced (since each $Y\cap X^{wP})$ is), and
Cohen-Macaulay (since $\pi_Y$ is flat and $\cY$ is Cohen-Macaulay). 
Thus, $\pi^{-1}_Y(0,\ldots,0)$ is reduced. By semicontinuity, it
follows that all fibers are reduced.
\end{proof}

For $Y$ satisfying the assumptions of Theorem \ref{degeneration}, let
$$
Y_0=\pi_Y^{-1}(0,\ldots,0)=\bigcup_{w\in W^P} X_{wP}\times (Y\cap X^{wP})
$$
regarded as a subvariety of $X\times X$.

\begin{corollary}\label{filtrations}
(i) With preceding notation and assumptions, the structure sheaf
$\cO_{Y_0}$ has an ascending filtration with associated graded
$$
\bigoplus_{w\in W^P} \cO_{X_{wP}}\otimes_{\mC} 
\cO_{Y\cap X^{wP}}(-Y\cap\partial X^{wP}).
$$
It has also a descending filtration with associated graded
$$
\bigoplus_{w\in W^P} \cO_{X_{wP}}(-\partial X_{wP})\otimes_{\mC} 
\cO_{Y\cap X^{wP}}.
$$

\noindent
(ii) If in addition $Y$ is normal, then the dualizing
sheaf $\omega_{Y_0}$ has a descending filtration with associated
graded 
$$
\bigoplus_{w\in W^P} \omega_{X_{wP}}\otimes_{\mC} 
\omega_{Y\cap X^{wP}}(Y\cap \partial X^{wP}).
$$
It has also an ascending filtration with associated graded 
$$
\bigoplus_{w\in W^P} 
\omega_{X_{wP}}(\partial X_{wP})\otimes_{\mC}\omega_{Y\cap X^{wP}}. 
$$
\end{corollary}

\begin{proof}
(i) We adapt the argument of \cite{BP} Theorem 11 to this setting. We
may index the finite poset $W^P=\{w_1,\ldots,w_N\}$ so that $i\leq j$
whenever $w_i\leq w_j$. Let 
$$
Z_i=X_{w_iP}\times (Y\cap X^{w_iP}),~
Z_{\geq i}=\cup_{j\geq i}\, Z_j \text{ and } 
Z_{>i}=\cup_{j>i}\, Z_j
$$ 
for $1\leq i\leq N$. Then $Z_{\geq 1}=Y_0$ and 
$Z_{\geq N}=X_{w_N}\times (Y\cap X^{w_N})$. We claim that
$$
Z_i\cap Z_{>i}= X_{w_iP}\times (Y\cap\partial X^{w_iP}).
$$
To see this, we may assume that $Y=X$. Then $Z_i\cap Z_{>i}$ is a
union of products $X_{uP}\times X^{vP}$ for certain $u,v$ in $W^P$. We must
have $u\leq w_i\leq v$ (since $X_{uP}\times X^{vP}\subseteq Z_i$) and
$w_i\neq v$ (since $X_{uP}\times X^{vP}\subseteq Z_{>i}$). Thus, 
$Z_i\cap Z_{>i}$ is contained in $X_{w_iP}\times \partial X^{w_iP}$. 
Conversely, if $X^{vP}\subseteq \partial X^{w_iP}$, then $v=w_j$ for some
$j$ such that $w_j>w_i$, whence $j>i$ ; this yields the opposite
inclusion. The claim is proved. 

Now consider the exact sequence
$$
0\to\cI_i\to\cO_{Z_{\geq i}}\to\cO_{Z_{>i}}\to 0
$$
where $\cI_i$ denotes the ideal sheaf of $Z_{>i}$ in 
$Z_{\geq i}$. Then $\cI_i$ identifies with the ideal sheaf of 
$Z_i\cap Z_{>i}$ in $Z_i$ ; by the claim, this is the ideal sheaf of 
$X_{w_iP}\times (Y\cap \partial X^{w_iP})$ in 
$X_{w_iP}\times (Y\cap X^{w_iP})$. This yields the ascending filtration
of $\cO_{Y_0}$.

With obvious notation, we obtain likewise
$$
Z_i\cap Z_{<i}=\partial X_{w_iP}\times (Y\cap X^{w_iP})
$$
which yields the descending filtration.

\noindent
(ii) By Lemma \ref{boundary}, the sheaf $\cI_i$ is Cohen-Macaulay of
depth $\dim(Y)$. Now a descending induction on $i$ shows that each 
$Z_{\geq i}$ is a Cohen-Macaulay variety of dimension
$\dim(Y)$. Furthermore, we obtain exact sequences
$$\displaylines{
0\to 
Ext^{\dim(X\times X)-\dim(Y)}_{X\times X}
(\cO_{Z_{>i}},\omega_{X\times X})
\to
Ext^{\dim(X\times X)-\dim(Y)}_{X\times X}
(\cO_{Z_{\geq i}},\omega_{X\times X})\to
\hfill\cr\hfill
\to 
Ext^{\dim(X\times X)-\dim(Y)}_{X\times X}(\cI_i,\omega_{X\times X})
\to 0,
\cr}$$
that is, 
$$
0\to \omega_{Z_{>i}} \to \omega_{Z_{\geq i}} \to
\omega_{X_{w_iP}}\otimes_{\mC} Ext_X^{\codim(Y\cap X^{w_iP})}
(\cO_{Y\cap X^{w_iP}}(-Y\cap\partial X^{w_iP}),\omega_X)
\to 0.
$$
Since $Y\cap X^{w_iP}$ and $Y\cap \partial X^{w_iP}$ are Cohen-Macaulay
by Lemma \ref{boundary}, we obtain
$$\displaylines{
Ext_X^{\codim(Y\cap X^{w_iP})}
(\cO_{Y\cap X^{w_iP}}(-Y\cap\partial X^{w_iP}),\omega_X)
\hfill\cr\hfill
\cong
Hom_X(\cO_{Y\cap X^{w_iP}}(-Y\cap\partial X^{w_iP}),
\omega_{Y\cap X^{w_iP}}).
\cr}$$
The latter is isomorphic to 
$\omega_{Y\cap X^{w_iP}}(Y\cap\partial X^{w_iP})$, by Lemma
\ref{boundary} again. This yields a descending filtration of 
$\omega_{Z_{\geq 1}}=\omega_{Y_0}$, with associated graded as
claimed. The ascending filtration is obtained by replacing 
$Z_{\geq i}$ with $Z_{\leq i}$.
\end{proof}

Next we derive from Corollary \ref{filtrations} several formulae for
decomposing $[\cO_Y]$ and $[\omega_Y]$ in the Grothendieck group
$K(X)$. Recall that this group is freely generated by the classes
$[\cO_{X_{wP}}]$ where $w\in W^P$ (see \cite{KK} \S 4). Another natural
basis of $K(X)$ consists of the classes
$[\cO_{X_{wP}}(-\partial X_{wP})]=[\cO_{X_{wP}}]-[\cO_{\partial X_{wP}}]$.
Using the duality involution $*$, we obtain two additional bases: 
the $[\omega_{X_{wP}}]$, and the $[\omega_{X_{wP}}(\partial X_{wP})]$.

\begin{corollary}\label{formulae}
For any Cohen-Macaulay closed subvariety $Y$ of $X$, we have in $K(X)$:
$$
[\cO_Y]=\sum_{w\in W^P} 
\chi(\cO_{Y\cap g X^{wP}}(-Y\cap g\partial X^{wP}))[\cO_{X_{wP}}]
=\sum_{w\in W^P} 
\chi(\cO_{Y\cap g X^{wP}})[\cO_{X_{wP}}(-\partial X_{wP})]
$$
for general $g\in G$. If in addition $Y$ is normal, then
$$
[\omega_Y]=\sum_{w\in W^P} 
\chi(\omega_{Y\cap gX^{wP}}(Y\cap g\partial X^{wP}))
\,[\omega_{X_{wP}}]
=\sum_{w\in W^P} 
\chi(\omega_{Y\cap g X^{wP}})\,[\omega_{X_{wP}}(\partial X_{wP})].
$$
\end{corollary}

\begin{proof}
Since $\diag(Y)$ and $Y_0$ are two fibers of the flat family $\cY$
over the affine space $\mA^r$, we have $[\cO_{\diag(Y)}]=[\cO_{Y_0}]$
in $K(X\times X)$. By Corollary \ref{filtrations} (i), it follows that
$$\displaylines{
[\cO_{\diag(Y)}]=\sum_{w\in W^P} [\cO_{X_{wP}}\otimes_{\mC}
\cO_{Y\cap gX^{wP}}(-Y\cap\partial gX^{wP})]
\hfill\cr\hfill
=\sum_{w\in W^P} [\cO_{X_{wP}}(-\partial gX_{wP})\otimes_{\mC} 
\cO_{Y\cap gX^{wP}}].
\cr}$$
Now let 
$$
p_1:X\times X\to X
$$ 
be the projection to the first factor, and denote 
$$
p_{1*}:K(X\times X)\to K(X)
$$ 
the corresponding pushforward map. Then
$p_{1*}[\cO_{\diag(Y)}]=[\cO_Y]$, whereas 
$$
p_{1*}([\cF\otimes_{\mC}\cG]) = \chi(\cG)\,[\cF]
$$
for all coherent sheaves $\cF$, $\cG$ on $X$. This yields our formulae
for $\cO_Y$.

To obtain the formulae for $\omega_Y$, notice that the sheaf
$\omega_{\cY}$ is flat over $\mA^r$, since $\cY$ is Cohen-Macaulay and
flat over $\mA^r$. Furthermore, the restriction to $\omega_{\cY}$ to
the fiber at $(1,\ldots,1)$ (resp. $(0,\ldots,0)$) is isomorphic to
$\omega_{\diag(Y)}$ (resp. $\omega_{Y_0}$). Thus,
$\omega_{\diag(Y)}=\omega_{Y_0}$ in $K(X\times X)$. Now both
formulae follow from Corollary \ref{filtrations} (ii) by the
preceding argument. 

Alternatively, these formulae may be derived from those for $[\cO_Y]$
by applying the involution $*$ and duality in $Y\cap gX^{wP}$. For, using
Lemma \ref{boundary}, we obtain isomorphisms
$$
H^i(\cO_{Y\cap gX^{wP}}(-Y\cap g\partial X^{wP}))^* 
= H^{\dim(Y\cap gX^{wP})-i}(\omega_{Y\cap gX^{wP}}(Y\cap g\partial X^{wP})),
$$
whence 
$$
\chi(\cO_{Y\cap gX^{wP}}(-Y\cap g\partial X^{wP}))=
(-1)^{\dim(Y\cap gX^{wP})}
\chi(\omega_{Y\cap gX^{wP}}(Y\cap g\partial X^{wP})).
$$
\end{proof}

\section{A vanishing theorem}

Consider a normal Cohen-Macaulay subvariety $Y$ of a flag variety $X$.  
In Section 2, we constructed a flat family in $X\times X$ with general
fibers certain translates of $\diag(Y)$. Furthermore, the special fiber
$Y_0$ is Cohen-Macaulay, and its canonical sheaf has a filtration with
layers the
$\omega_{X_{wP}}\otimes_{\mC}
\omega_{Y\cap gX^{wP}}(Y\cap g\partial X^{wP})$,
where $g$ is a general element of $G$. 

Let $p_1:X\times X\to X$ be the first projection. If all higher
cohomology groups of the sheaves 
$\omega_{Y\cap gX^{wP}}(Y\cap g\partial X^{wP})$ vanish, then 
$R^ip_{1*}(\omega_{Y_0})=0$ for $i\geq 1$. Furthermore, the sheaf
$p_{1*}\omega_{Y_0}$ has a filtration with layers the
$\omega_{X_{wP}}$ of respective multiplicities 
$h^0(\omega_{Y\cap gX^{wP}}(Y\cap g\partial X^{wP}))$. 
Thus, the following equalities hold in $K(X)$:
$$\displaylines{
p_{1*}[\omega_{Y_0}]=\sum_{w\in W^P}\,
h^0(\omega_{Y\cap gX^{wP}}(Y\cap g\partial X^{wP}))
\, [\omega_{X_{wP}}],
\hfill\cr\hfill
p_{1*}[\cO_{Y_0}]=\sum_{w\in W^P} (-1)^{\codim(X_{wP})-\codim(Y)}
h^0(\omega_{Y\cap gX^{wP}}(Y\cap g\partial X^{wP}))
\, [\cO_{X_{wP}}].
\cr}$$
On the other hand, since $p_1$ maps isomorphically $\diag(Y)$ to $Y$,
we obtain $p_{1*}[\cO_{Y_0}]=[\cO_Y]$ and
$p_{1*}[\omega_{Y_0}]=[\omega_Y]$ ; together with the preceding
equalities, this yields a sharper version of Theorem \ref{signs}. 

By our next result, this vanishing condition holds in the case
where $Y$ has rational singularities. In fact we shall prove a
slightly more general vanishing theorem, in view of further
applications in the next section.

\begin{theorem}\label{vanishing}
Consider a flag variety $X=G/P$, a closed subvariety $Y$ and an
invertible sheaf $\cL=\cL_X(\lambda)$. Assume that $Y$ has rational
singularities, and that $\lambda$ is dominant. Then we have for 
general $g\in G$:
$$
H^i(\cL_{Y\cap g X^{wP}}(-\lambda)(-Y\cap g\partial X^{wP}))=0
$$
for all $w\in W^P$ and $i<\dim(Y\cap g X^{wP})=\codim(X_{wP})-\codim(Y)$. 
Equivalently,
$$
H^i(\cL_{Y\cap g X^{wP}}(\lambda)\otimes
\omega_{Y\cap gX^{wP}}(Y\cap g\partial X^{wP}))=0
$$
for all $w\in W^P$ and $i\geq 1$.
\end{theorem}

\begin{proof} We first consider the case where $X$ is the full flag
variety ; furthermore, we replace $X^w$ by $X_w$ for simplicity. 

Recall that each intersection $Y\cap g X_w$ is the fiber
$\pi^{-1}(g)$, with notation displayed by the commutative diagram
$$
\CD
 G @<{\pi}<< Z @>{\mu}>> Y\\
 @V{id}VV    @V{\iota}VV @V{i}VV\\
 G @<{p}<<   G\times X_w @>{m}>> X
\endCD
$$
where the square on the right is cartesian. Recall also that
$Z$ has rational singularities. Let 
$$
\partial Z=Y\times _X (G\times\partial X_w)
$$
a subvariety of codimension $1$ in $Z=Y\times_X (G\times X_w)$. 
For general $g\in G$, we have
$$
(\mu^*\cL_Y(\lambda)\otimes\omega_Z(\partial Z))\vert_{\pi^{-1}(g)}
\cong
\cL_{Y\cap gX_w}(\lambda)\otimes
\omega_{Y\cap gX_w}(Y\cap g\partial X_w).
$$
Thus, our statement (ii) is a consequence of the following assertion:
$$
R^i\pi_*(\mu^*\cL_Y(\lambda)\otimes\omega_Z(\partial Z))=0 
\text{ for }i\geq 1.\eqno(1)
$$

We shall deduce (1) from the Kawamata-Viehweg vanishing theorem,
like in the proof of Lemma \ref{boundary}. Since that theorem concerns
nonsingular varieties, we first construct a desingularization of $Z$.

Let $\varphi:\tY\to Y$ be a desingularization. On the other
hand, let $\psi:\tX_w\to X_w$ be a standard desingularization as
in the proof of Lemma \ref{boundary}. 
Composing $id\times\psi: G\times \tX_w\to G\times X_w$ with 
the multiplication map $m:G\times X_w\to X$ defines 
$\tilde m:G\times\tX_w\to X$. Define likewise 
$\tilde i=i\circ \varphi$, $\tilde p=p\circ(id\times\psi)$ and
consider the commutative diagram
$$
\CD
 G @<{\tilde\pi}<< \tZ @>{\tilde\mu}>> \tY\\
 @V{id}VV            @V{\tilde\iota}VV @V{\tilde i}VV\\
 G @<{\tilde p}<<   G\times \tX_w @>{\tilde m}>> X
\endCD
$$
where the square on the right is cartesian. Since $\tilde m$ is a
$G$-equivariant morphism from the nonsingular variety 
$G\times\tX_w$ to $X=G/B$, it is a locally trivial fibration with
nonsingular fiber. Thus, the same holds for $\tilde\mu$, so that
$\tZ$ is nonsingular as well. The map 
$id\times\tilde f\times\varphi:G\times\tX_w\times\tY\to
G\times X_w\times Y$ 
is a desingularization ; it restricts to a proper morphism 
$$
f:\tZ\to Z
$$
which is clearly birational. Thus, $f$ is a desingularization of $Z$.
The subset
$$
\partial\tZ=\tY\times_X (G\times\partial\tX_w)
$$
is a union of nonsingular irreducible divisors intersecting
transversally in $\tZ$ ; clearly, $f(\partial\tZ)=\partial Z$.

We claim that
$$
\omega_{\tZ}(\partial\tZ)=
\tilde\mu^*(\omega_{\tY}\otimes\varphi^*\cL_Y(\rho)).
$$
To verify this, notice that
$$
\omega_{\tZ}=\tilde\mu^*\omega_{\tY}\otimes\omega_{\tZ/\tY},
$$
since $\mu$ is a locally trivial fibration. Furthermore, 
$$
\omega_{\tZ/\tY}=\tilde\iota^*\omega_{(G\times \tX_w)/X}
=\tilde\iota^*(\omega_{G\times \tX_w}\otimes \tilde m^*\omega_X^{-1})
=\tilde\iota^*(\cO_{G\times \tX_w}(-G\times\partial \tX_w)
\otimes \tilde m^*\cL_X(\rho)),
$$
since $\omega_G=\cO_G$ and 
$\omega_{\tX_w}=\cO_{\tX_w}(-\partial\tX_w)\otimes
{\tilde m}^*\cL_X(-\rho)$. Therefore, 
$$
\omega_{\tZ}=\tilde\mu^*\omega_{\tY}\otimes\cO_{\tZ}(-\partial\tZ)
\otimes\tilde\iota^*\tilde m^*\cL_X(\rho)
$$
which implies the claim.

We next obtain the analogue of (1) for $\tZ$, that is,
$$
R^i\tilde\pi_*(\tilde\mu^*\varphi^*\cL_Y(\lambda)\otimes
\omega_{\tZ}(\partial\tZ))=0 \text{ for }i\geq 1.
\eqno(2)
$$
This follows from a relative version of the Kawamata-Viehweg vanishing
theorem. Specifically, recall that $\partial\tX_w$ is the support of
a very ample divisor  $b_1\tX_1+\cdots+b_{\ell}\tX_{\ell}$ of $\tX_w$,
with the notation of the proof of Lemma \ref{boundary}. Let
$\tZ_i=\tY\times_X(G\times\tX_i)$ for $i=1,\ldots,\ell$ ; then the
nonempty $\tZ_i$ are the irreducible components of $\partial\tZ$. 
Define positive integers $N,a_1,\ldots,a_{\ell}$ as in the proof of
Lemma \ref{boundary} and let 
$\cM=(\tilde\mu^*\varphi^*\cL_Y(\lambda))(\partial\tZ)$ and
$D=a_1\tZ_1+\cdots+a_{\ell}\tZ_{\ell}$. 
Then the invertible sheaf
$$
\cM^N(-D)=(\tilde\mu^*\tilde i^*\cL_X(N\lambda))
(b_1\tZ_1+\cdots+b_{\ell}\tZ_{\ell})
=\tilde\iota^*(\tilde m^*\cL_X(N\lambda)
(G\times(b_1\tX_1+\cdots+b_{\ell}\tX_{\ell})))
$$
is the pullback under $\tilde\iota$ of a very ample invertible sheaf
on $G\times\tX_w$. Since $\tilde\iota$ is generically injective, it
follows that $\cM^N(-D)$ is $\tilde\pi$-numerically effective and
$\tilde\pi$-big. Therefore, $R^i\tilde\pi_*(\cM\otimes\omega_{\tZ})=0$
for $i\geq 1$, by \cite{EV} Corollary 6.11. This proves (2).

Likewise, $\cM^N(-D)$ is $f$-numerically effective, so that 
$$
R^i f_*(\tilde\mu^*\varphi^*\cL_Y(\lambda)
\otimes\omega_{\tZ}(\partial\tZ))=0 \text{ for }i\geq 1.\eqno(3)
$$
Finally, we claim that
$$
f_*(\tilde\mu^*\varphi^*\cL_Y(\lambda)
\otimes\omega_{\tZ}(\partial \tZ))
=\mu^*\cL_Y(\lambda)\otimes\omega_Z(\partial Z).\eqno(4)
$$
Together with (2), (3) and the Leray spectral sequence for
$\tilde\pi=\pi\circ f$, this will imply assertion (1). 

To check (4), we factor $f$ into $\varphi'\circ f'$, with notation
displayed in the commutative diagram
$$\CD
 \tZ @>{f'}>> Z' @>{\varphi'}>> Z\\
 @V{\tilde\mu}VV @V{\mu'}VV   @V{\mu}VV\\
 \tY @>{id}>>   \tY  @>{\varphi}>> Y
\endCD$$
where the square on the right is cartesian. Notice that $Z'$ has
rational singularities, and that $f'$ is a desingularization ;
furthermore, we obtain
$$
\omega_{Z'}(\partial Z')=
\mu^{'*}(\omega_{\tY}\otimes\varphi^*\cL_Y(\rho)),
$$
by the preceding arguments for determining $\omega_{\tZ}(\partial\tZ)$, 
applied to the regular locus of $Z'$. Thus,
$$
f'_*(\tilde\mu^*\varphi^*\cL_Y(\lambda)\otimes\omega_{\tZ}(\partial\tZ))=
f'_*f^{'*}\mu^{'*}(\omega_{\tY}\otimes\varphi^*\cL_Y(\lambda+\rho))
=\mu^{'*}(\omega_{\tY}\otimes\varphi^*\cL_Y(\lambda+\rho)).
$$
It follows that
$$
f_*(\tilde\mu^*\varphi^*\cL_Y(\lambda)
\otimes\omega_{\tZ}(\partial \tZ))
=
\varphi'_*\mu^{'*}(\omega_{\tY}\otimes\varphi^*\cL_Y(\lambda+\rho))=
\mu^*\varphi_*(\omega_{\tY}\otimes\varphi^*\cL_Y(\lambda+\rho)),
$$
where the latter equality holds since $\mu$ is flat. By the projection
formula and rationality of singularities of $Y$, this yields
$$
f_*(\tilde\mu^*\varphi^*\cL_Y(\lambda)
\otimes\omega_{\tZ}(\partial \tZ))=
\mu^*(\omega_Y\otimes\cL_Y(\lambda+\rho)).
$$
And one may check as above that the latter equals
$\mu^*\cL_Y(\lambda)\otimes\omega_Z(\partial Z)$.
This completes the proof of (4) and hence of Theorem 
\ref{vanishing}, in the case where $X=G/B$.

Finally, in the case where $X=G/P$, one argues by reducing to $G/B$ as
in the proof of Lemma \ref{boundary} ; we skip the details.
\end{proof}

\medskip

\noindent
{\sl Remark.} Consider an arbitrary closed subvariety $Y$ of $X$ and a
desingularization $\varphi:\tY\to Y$. Then the sheaves
$R^i\varphi_*\cO_{\tY}$ are independent of the choice of $\varphi$, so
that the same holds for the class 
$$
\varphi_*[\cO_{\tY}]=\sum_{i=0}^{\dim(Y)}\,
(-1)^i\,[R^i\varphi_*\cO_{\tY}]
$$ 
in $K(X)$. Thus, we may define integers $b_Y^w$ by
$$
\varphi_*[\cO_{\tY}]=\sum_{w\in W^P}\,b_Y^w\,[\cO_{X_{wP}}].
$$
Then one may adapt the arguments of Sections 1 and 2 to obtain
$$
b_Y^w=\chi(\cO_{\tY\times_X gX^{wP}}(-\tY\times_X g\partial X^{wP}))
$$
for general $g\in G$. And the proof of Theorem \ref{vanishing}
actually shows that
$$
H^i(\cO_{\tY\times_X gX^{wP}}(-\tY\times_X g\partial X^{wP}))=0
~\text{ for }~i<\codim(X_{wP})-\codim(Y). 
$$
As a consequence,
$$
(-1)^{\codim(X_{wP})-\codim(Y)}\, b_Y^w\geq 0
$$
for all $w\in W^P$.

This admits a simpler formulation in terms of the sheaf
$\varphi_*\omega_{\tY}$. The latter is also independent of the choice 
of $\varphi$, and is called the sheaf of absolutely regular
differential forms on $Y$ ; we denote it by $\tilde\omega_Y$. 
Furthermore, $R^i\varphi_*\omega_{\tY}=0$ for $i\geq 1$ (see
\cite{GR}). Using duality for the morphism $\tY\to X$, it follows that
$$
[\tilde\omega_Y]=\sum_{w\in W^P}\,
(-1)^{\codim(X_{wP})-\codim(Y)}\, b_Y^w\,[\omega_{X_w}].
$$
In other words, the class of $\tilde\omega_Y$ is a nonnegative
combination of classes of dualizing sheaves of Schubert varieties. 
This generalizes Theorem \ref{signs} to all closed subvarieties of
flag varieties.

\section{Restricting homogeneous line bundles to Schubert subvarieties}

The Grothendieck ring of the full flag variety $X$ is generated as an
additive group by classes of invertible sheaves, see \cite{Mar}. This
raises the question of describing the product of such classes with
classes of structure sheaves of Schubert varieties. For any weight
$\lambda$ and for any $v$ in $W$, we have in $K(X)$:
$$
[\cL_X(\lambda)] \cdot [\cO_{X_v}]=[\cL_{X_v}(\lambda)]
=\sum_{w\in W} \, c_v^w(\lambda)\, [\cO_{X_w}]
$$
for uniquely defined integer coefficients $c_w^v(\lambda)$. Our next
result expresses these coefficients in geometric terms.

\begin{theorem}\label{line}
For any weight $\lambda$ and for any $v\in W$, we have
$$
c_v^w(\lambda)=
\chi(\cL_{X_v\cap X^w}(\lambda)(-X_v\cap\partial X^w)).
$$
As a consequence, we have the duality formula
$$
c_v^w(-\lambda)=c_{w_ow}^{w_ov}(-w_o\lambda).
$$
If in addition $\lambda$ is dominant, then
$$
c_v^w(\lambda)=
h^0(\cL_{X_v\cap X^w}(\lambda)(-X_v\cap \partial X^w)).
$$
\end{theorem}

\begin{proof}
We apply Theorem \ref{degeneration} and Corollary \ref{filtrations} to
the normal, Cohen-Macaulay variety $X_v$. This yields
$$
[\cO_{\diag(X_v)}]=\sum_{w\in W}\,[\cO_{X_w}]\times
[\cO_{X_v\cap X^w}(-X_v\cap \partial X^w)]
$$
in $K(X\times X)$. Multiplying both sides by $p_2^*[\cL_X[\lambda)]$
(where $p_1,p_2:X\times X\to X$ are the projections) and then applying
$p_{1*}$, we obtain 
$$
p_{1*}([\cO_{\diag(X_v)}]\cdot p_2^*[\cL_X(\lambda)])=
\sum_{w\in W} \, 
\chi(\cL_{X_v\cap X^w}(\lambda)(-X_v\cap \partial X^w))\,[\cO_{X_w}]
$$
for any weight $\lambda$. But since $p_1:\diag(X_v)\to X_v$ is an
isomorphism, we have
$$
p_{1*}([\cO_{\diag(X_v)}]\cdot p_2^*[\cL_X(\lambda)])=
p_{1*}[\diag(\cL_{X_v}(\lambda))]=[\cL_{X_v}(\lambda)].
$$
This proves our first formula.

Recalling that
$$
\omega_{X_v\cap X^w}=\cO_{X_v\cap X^w}
(-(X_v\cap \partial X^w)-(\partial X_v\cap X^w))
$$
as seen in Section 1, we obtain
$$
\chi(\cL_{X_v\cap X^w}(\lambda)(-X_v\cap \partial X^w))=
\chi(\cL_{X_v\cap X^w}(\lambda)\otimes
\omega_{X_v\cap X^w}(\partial X_v\cap X^w)).
$$
By Lemma \ref{boundary} and duality in the variety $X_v\cap X^w$ of
dimension $\ell(v)-\ell(w)$, it follows that 
$$
c_v^w(\lambda)=(-1)^{\ell(v)-\ell(w)}\,
\chi(\cL_{X_v\cap X^w}(-\lambda)(-\partial X_v\cap X^w)).
$$
Since $X_v=w_oX^{w_ov}$ and $X^w=w_oX_{w_ow}$, this implies our second
formula.

If in addition $\lambda$ is dominant, then
$H^i(\cL_{X_v\cap X^w}(\lambda)\otimes
\omega_{X_v\cap X^w}(\partial X_v\cap X^w))=0$ for every $i\geq 1$,
as follows from Theorem \ref{vanishing}. This yields our third
formula.
\end{proof}

As a consequence, $c_v^w(\lambda)=0$ unless $w\leq v$, and
$c_v^v(\lambda)=1$ ; furthermore, $c_v^w(\lambda)\geq 0$ if $\lambda$
is dominant. 

The definition of the coefficients $c_v^w(\lambda)$ implies that 
$$
c_v^w(\lambda+\mu)=\sum_{x\in W,\,w\leq x\leq v}\, 
c_v^x(\lambda) \,c_x^w(\mu)
$$
for all weights $\lambda$ and $\mu$. Together with the second formula
in Theorem \ref{line}, this shows that the $c_v^w(\lambda)$ may be
expressed in terms of the $c_v^w(\omega_i)$, where the $\omega_i$ are
the fundamental weights. The latter are related to certain structure
constants $c_{u,v}^w$ as follows.

\begin{lemma}
With preceding notation, we have for $v\neq w$:
$$
c_v^w(-\omega_i)=-c_{w_os_i,v}^w\text{ and }
c_v^w(\omega_i)=(-1)^{\ell(v)-\ell(w)-1}\,
c_{s_iw_o,w_ov}^{w_ow}.
$$
\end{lemma}

\begin{proof}
The invertible sheaf $\cL_X(\omega_i)$ has a section with zero
subscheme the Schubert variety $X_{w_os_i}$. This yields an exact
sequence 
$$
0\to \cL_X(-\omega_i) \to \cO_X \to \cO_{X_{w_os_i}} \to 0,
$$
whence  $[\cL_X(-\omega_i)]=[\cO_X] - [\cO_{X_{w_os_i}}]$ in $K(X)$. 
Multiplying this equality by $[\cO_{X_v}]$ yields
$$
[\cL_{X_v}(-\omega_i)]=
[\cO_{X_v}]-\sum_{w\in W} \, c_{w_os_i,v}^w\,[\cO_{X_w}]
$$
which implies our first formula. The second formula follows by duality.
\end{proof}

Notice that all results of this section extend to the setting of
$T$-equivariant $K$-theory (see \cite{KK}). For the intersections
$X_v\cap X^w$ are invariant under the action of $T$ on $X$, and the
constructions of Section 2 are equivariant with respect to this
action. This yields a geometric proof for the positivity of
$[\cL_{X_v}(\lambda)]$ in $K^T(X)$, due to Pittie and Ram (see 
\cite{PR} Corollary, p. 106) and Mathieu \cite{Mat}. And this raises
the question of a positivity result for the structure constants of
$K^T(X)$ ; see \cite{G} for such a result in the setting of
$T$-equivariant cohomology.

\end{document}